\input amstex
\documentstyle{amsppt}
\TagsOnRight
\def\N{\Bbb N}
\def\Z{\Bbb Z}

\def\prob{\bold P}
\def\expect{\bold E}

\def\ind{1\!\!1}
\def\wt{\widetilde}
\def\vareps{\varepsilon}


\topmatter
\title
No More than Three Favourite Sites for Simple Random Walk
\endtitle

\leftheadtext{B\'alint T\'oth}
\rightheadtext{No more than three favourite sites for simple random walk}
\author
B\'alint T\'oth
\endauthor
\affil
Technical University Budapest, 
Institute of Mathematics
\\
Egry J. u. 1, \ \ 
H-1111 Budapest, \ \ 
Hungary
\\
e-mail: balint\@math.bme.hu
\endaffil
\abstract
{
We prove that, with probability one, eventually  there are no
more than three favourite 
(i.e. most visited)
sites of simple symmetric random walk.
This partially answers a relatively 
long standing question of 
P\'al Erd\H os and P\'al R\'ev\'esz. 
\newline\newline
{\smc Key Words:} 
Random walks, 
local time, 
favourite sites,  
most visited sites
\newline
{\smc 1991 AMS Subject Classification:}  
60J15, 60J55}
\endabstract
\endtopmatter


\document
\advance\baselineskip by 3pt
\advance\parindent by 3mm



\bigskip
\noindent
{\bf 1. Introduction and Main Result.}
 
\medskip
\noindent
Let $S(t)$, $t\in\Z_+$ be  a simple symmetric random walk on
$\Z$ with initial state $S(0)=0$. Its {\it upcrossings}, {\it
downcrossings} and {\it (site) local time} are defined for
$t\in\N$ and $x\in\Z$ as follows:  
$$
\align
U(t,x)
&:=
\#\{0<s\le t: S(s)=x, \ S(s-1)=x-1\},
\tag 1.1
\\
D(t,x)
&:=
\#\{0<s\le t: S(s)=x, \ S(s-1)=x+1\},
\tag 1.2
\\
L(t,x)
&:=
\#\{0<s\le t: S(s)=x\}
=
U(t,x)+D(t,x).
\tag 1.3
\endalign
$$
The following identities are  straightforward:
$$
\align
U(t,x)-D(t,x-1) 
&=
\ind_{\{0<x\le S(t)\}} - \ind_{\{S(t)<x\le0\}}
\tag 1.4
\\
D(t,x)-U(t,x+1)
&=
\ind_{\{S(t)\le x <0\}} - 
\ind_{\{0\le x < S(t)\}}.
\tag 1.5
\endalign
$$
And from these it follows that
$$
\align
L(t,x)
&=
D(t,x)+D(t,x-1)+
\ind_{\{0<x\le S(t)\}} - \ind_{\{S(t)<x\le 0\}}
\\
&=
U(t,x)+U(t,x+1)+
\ind_{\{S(t)\le x < 0\}} -
\ind_{\{0\le x<  S(t)\}}.
\tag 1.6
\endalign
$$
The set of {\it favourite (or: most visited) sites\/} of the random walk at time
$t\in\N$, are those sites where the local time attains its
maximum value:
$$
\Cal K(t) :=
\{y\in\Z: L(t,y)=\max_{z\in\Z}L(t,z)\}.
\tag 1.7
$$
It is clear that the number of favourite sites changes in time
as follows:
$$
\#\Cal K(t+1)
=
\left\{\matrix\format\l\quad&\quad\c\quad&\quad\l\\
\#\Cal K(t)   
& 
\text{ if } 
& 
S(t+1)\notin\Cal K(t+1)   
\\
\#\Cal K(t)+1 
& 
\text{ if } 
& 
\Cal K(t+1)=\Cal K(t)\cup\{S(t+1)\} 
\\ 
1     
& 
\text{ if } 
& 
\Cal K(t+1)=\{S(t+1)\}\subset\Cal K(t).
\endmatrix\right.
\tag 1.8
$$
In plain words one of the following three possibilities can
occur at each step of the walk: Either the currently occupied
site is not favourite and $\Cal K$ remains unchanged. Or the
currently occupied site becomes a new favourite {\it beside\/}
the favourites of the previous stage, and thus the number of
favourites increases by one. Or, finally, a favourite site is
visited and so this site becomes now the only new favourite. No
other possibility. 

Clearly $\#\Cal K(t)\ge1$ for all $t\ge 1$, and it is easy to
verify that for infinitely many times, $t\ge1$, there are at least
two favourite sites: $\#\Cal K(t)\ge2$.  P\'al Erd\H os and
P\'al R\'ev\'esz formulated and repeatedly raised the following

\medskip
\noindent
{\bf Question:}
Does it happen that $\#\Cal K(t)\ge r$ infinitely often (i.e.
almost surely for infinitely many times $t\ge1$) for $r=3,4,
\dots$? 

\medskip
\noindent
See e.g. Erd\H os and R\'ev\'esz (1984), (1987), (1991), Erd\H os (1994) or  
R\'ev\'esz (1990) for an extended list of related questions and 
problems. 

Questions related to the asymptotic behaviour of the favourite 
(or most visited) sites of a random walk were considered by many 
authors since the mid-eighties. We quote here a few relevant results, 
with no claim of exhaustiveness. 
\newline
$\bullet$
Bass and Griffin (1985) prove that almost surely, the set of 
favourites is transient. 
More exactly: they prove that the distance of the set of favourite 
sites from the origin increases faster than $\sqrt{n}/(\log n)^{11}$ 
but slower than $\sqrt{n}/(\log n)$. 
\newline
$\bullet$
Cs\'aki and Shi (1998) prove that the distance between the edge of 
the range of the random walk and the set of favourite sites 
increases as fast as $\sqrt{n}/(\log\log n)^{3/2}$
\newline
$\bullet$
Cs\'aki, R\'ev\'esz and Shi (2000) prove that the position 
of a favourite site can have as large jumps as 
$\sqrt{2 n \log\log n}$, i.e., comparable with the diameter of 
the full range of the random walk.  They also extend a much 
earlier result of Kesten (1965), identifying the set of joint 
limit points  of the set of favourite sites and the favourite 
values (i.e. max. values of local time), both rescaled by  
$\sqrt{2 n \log\log n}$
\newline
$\bullet$
There are many papers dealing with similar qestions in the context 
of symmetric stable processes rather than random walks 
(or Brownian motion). See, e.g., Eisenbaum (1997), Bass, Eisenbaum and Shi (2000) and other 
papers quoted there.

In the present paper we answer in the negative the 
question of Erd\H os and R\'ev\'esz  quoted above, 
for $r\ge4$: we prove that
with probability 1, there are at most finitely many times
$t\ge1$ when there are four or more favourite sites of the
random walk $S(t)$. In T\'oth and Werner (1997) a similar result 
was proved for the set of {\it favourite edges\/} rather than 
favourite sites.  The present paper deals with the original question 
of Erd\H os and R\'ev\'esz. The starting general ideas of the present 
paper (see Sections 1-3) are very close to those of T\'oth and Werner 
(1997). However: the details of the proof require more refined 
estimates and arguments. On the technical level (see Sections 4-6) 
this proof is rather different. 

For $r\ge 1$ denote by $f(r)$ the (possibly infinite) number of
steps, when the currently occupied site is one of the $r$ actual
favourites:
$$
f(r)
:=
\#\{t\ge1: 
[S(t)\in\Cal K(t)] 
\ \land \ 
[\#\Cal K(t)=r]
\}.
\tag 1.9
$$ 
From (1.8) it follows that for any $r\ge 1$, $f(r+1)\le f(r)$.
(Both sides of the inequality could be infinite.)

The main result of this paper is the following 

\proclaim{Theorem 1}
$$
\expect\big(f(4)\big)<\infty.
\tag 1.10
$$
\endproclaim

\noindent{\sl Remarks:} 
(1) 
From this theorem the negative answer to the question of Erd\H
os and R\'ev\'esz clearly follows, for the cases $r\ge4$. 
\newline
(2) 
The case $r=3$ remains open. From the proof of the above theorem
it becomes clear that $\expect\Big(f(3)\Big)=\infty$.
Nevertheless we conjecture that $f(3)<\infty$, almost surely.

The paper is organized as follows: In Section 2 we perform some 
straightforward manipulations (essentially: rearrangements of 
sums). In Section 3 we recall the Ray-Knight Theorems for the 
local times of simple random walks. In Section 4
first we express our relevant probabilities and expectations
(found in Section 2) 
in terms of the Galton-Watson processes arising with the 
Ray-Knight representation. Then we formulate Proposition 1,
stating
some upper bounds on these probabilities and expectations,
and using these bounds we prove Theorem 1. 
The proof of Proposition 1 is postponed to the end of Section 5.
In Section 5 four lemmas and, as their consequence, Proposition 1 
are proved. Throughout the technical parts of the proofs  
smaller, quite plausible statements are invoked. These are called 
Side-lemmas (1 to 7). Their proofs are postponed to Section 6. 

Throughout the paper, in various upper 
bounds, multiplicative constants, 
respectively, constants in exponential 
rates  will be denoted generically by 
$C$, respectively, by $\gamma$. The 
values of these constants may vary even 
within one proof, but we hope there is no 
danger of confusion. 


\bigskip
\noindent
{\bf 2. Preparations.}
 
\medskip
\noindent
In the following transformations the {\it inverse local
times\/}, defined below for $k\in\N$ and $x\in\Z$, will play an
essential r\^ole:  
$$
\align
T_{U}(k,x)
&:=
\inf\{t\ge1: U(t,x)=k\},
\tag 2.1
\\
T_{D}(k,x)
&:=
\inf\{t\ge1: D(t,x)=k\}.
\tag 2.2
\endalign
$$
It turns out that questions related to 
the local time are easier to handle if 
the random walk is observed at the random 
stopping times $T_{U}(k,x)$ and 
$T_{D}(k,x)$ rather than at deterministic 
times $t\ge 1$. This `combinatorial trick' 
has its origin in Knight (1963) and has 
been successfully applied in various 
contexts. See, e.g., K\"onig (1996), 
T\'oth (1995), (1996) and references cited there.
 
We express $f(4)$ with the help of some straightforward
rearrangements of summations:
$$
f(4)=\sum_{x\in\Z}\Big(u(x)+d(x)\Big)
\tag 2.3
$$
where
$$
\align
u(x)
:=&
\sum_{t=1}^\infty
\ind_{\{
S(t)=x, \ 
S(t-1)=x-1, \ 
x\in\Cal K(t), \ 
\#\Cal K(t)=4
\}}
\\
=&
\sum_{t=1}^\infty\sum_{k=1}^\infty
\ind_{\{
T_{U}(k,x)=t, \ 
x\in\Cal K(t), \ 
\#\Cal K(t)=4
\}}
\\
=&
\sum_{k=1}^\infty
\ind_{\{
x\in\Cal K(T_{U}(k,x)), \  
\#\Cal K(T_{U}(k,x))=
4\}}
\tag 2.4
\endalign
$$
and 
$$
\align
d(x)
:=&
\sum_{t=1}^\infty
\ind_{\{
S(t)=x, \ 
S(t-1)=x+1, \ 
x\in\Cal K(t), \ 
\#\Cal K(t)=4
\}}
\\
=&
\sum_{t=1}^\infty\sum_{k=1}^\infty
\ind_{\{
T_{D}(k,x)=t, \ 
x\in\Cal K(t), \ 
\#\Cal K(t)=4
\}}
\\
=&
\sum_{k=1}^\infty
\ind_{\{
x\in\Cal K(T_{D}(k,x)), \  
\#\Cal K(T_{D}(k,x))=4
\}}.
\tag 2.5
\endalign
$$
Clearly,
$$
u(x)
\buildrel{\text{law}}\over=
d(-x)
\tag 2.6
$$
and, consequently
$$
\expect\Big(f(4)\Big)
=
2\sum_{x=1}^\infty
\expect\Big(u(x)\Big)
+
2\sum_{x=0}^\infty
\expect\Big(d(x)\Big)
\tag 2.7
$$
with 
$$
\align
\expect\Big(u(x)\Big)
&=
\sum_{k=1}^\infty
\prob\Big(
[x\in\Cal K(T_{U}(k,x))]
 \ \land \   
[\#\Cal K(T_{U}(k,x))=4]
\Big)
\tag 2.8
\\
\expect\Big(d(x)\Big)
&=
\sum_{k=1}^\infty
\prob\Big(
[x\in\Cal K(T_{D}(k,x))]
 \ \land \   
[\#\Cal K(T_{D}(k,x))=4]
\Big)
\tag 2.9
\endalign
$$
We shall show in details that 
$$
\sum_{x=1}^\infty
\expect\Big(u(x)\Big)=
\sum_{x=1}^\infty
\sum_{k=1}^\infty
\prob\Big(
[x\in\Cal K(T_{U}(k,x))] 
\ \land \   
[\#\Cal K(T_{U}(k,x))=4]
\Big)
<\infty.
\tag 2.10
$$
The similar statement 
$$
\sum_{x=0}^\infty
\expect\Big(d(x)\Big)=
\sum_{x=0}^\infty
\sum_{k=1}^\infty
\prob\Big(
[x\in\Cal K(T_{D}(k,x))]
\ \land \   
[\#\Cal K(T_{D}(k,x))=4]
\Big)
<\infty
\tag 2.11
$$
can be  proved in an identical way.


\bigskip
\noindent
{\bf 3. Ray-Knight Representation.}
 
\medskip
\noindent
Throughout this paper we denote by $Y_t$ a  critical
branching process with geometric offspring distribution 
(Galton-Watson process) and
by $Z_t$ a critical branching process with geometric
offspring distribution and one intruder at each generation.
$Y_t$ and $Z_t$ are Markov chains with state space $\Z_+$ and
transition probabilities: 
$$
\align
\prob\Big( Y_{t+1}=j \Big| Y_{t}=i \Big)
=
\pi(i,j)
:=
&
\left\{\matrix\format\c\quad&\quad\c\quad&\c\\
\dsize
\delta_{0,j} 
& 
\text{ if } 
& 
i=0,
\\ \\
\dsize
2^{-i-j}\frac{(i+j-1)!}{(i-1)!j!} 
& 
\text{ if } 
& 
i>0.
\endmatrix\right. 
\tag 3.1
\\ \\
\prob\Big( Z_{t+1}=j \Big| Z_{t}=i \Big)
=
\rho(i,j)
:=
&\,\,\,\,
2^{-i-j-1}\frac{(i+j)!}{i!j!}.
\tag 3.2
\endalign
$$
Let $k\ge0$ and $x\ge1$ be fixed integers 
and define the following three
processes: 
\item{$\bullet$}
$Z_t$, $0\le t \le x-1$, is a Markov chain with transition
probabilities $\rho(i,j)$ and initial state $Z_0=k$; 
\item{$\bullet$}
$Y_t$, $-1\le t<\infty$, is a Markov chain with transition
probabilities  $\pi(i,j)$ and initial state $Y_{-1}=k$;
\item{$\bullet$}
Finally, $Y'_t$, $0\le t <\infty$, is another Markov chain with
the same transition probabilities $\pi(i,j)$ and initial state
$Y'_0=Z_{x-1}$. 

\noindent
The three chains are independent, except for the fact that $Y'$
starts from the terminal state of $Z$. Using these three
chains we patch together the process
$$
\Delta_{x,k}(y):=
\left\{
\matrix
\format \c \quad & \quad \c \quad & \r & \,\,\l\\
Z_{x-y-1}
&
\text{ if }
&
0 \le y \le 
&
x-1
\\
Y_{y-x}
&
\text{ if }
&
x-1 \le y \le  
&
\infty
\\
Y'_{-y}
&
\text{ if }
&
-\infty \le y \le 
&
0.
\endmatrix\right.
\tag 3.3
$$
We also define 
$$
\Lambda_{x,k}(y) :=
\Delta_{x,k}(y) + \Delta_{x,k}(y-1) + \ind_{\{0<y\le x\}}.
\tag 3.4
$$
According to the, by now classical, Ray-Knight Theorems on the
local time process of simple symmetric random walks on $\Z$ 
(cf. Knight (1963), Ray (1963)),
for any integers $x\ge1$ and $k\ge0$:
$$
\Big(
\Delta_{x,k}(y), \ \ y\in\Z
\Big) 
\buildrel{\text{law}}\over=
\Big( 
D(T_U(k+1,x),y), \ \ y\in\Z
\Big). 
\tag 3.5
$$
Using (1.6) and (3.4), from this we get 
$$
\Big
(\Lambda_{x,k}(y), \ \ y\in\Z
\Big) 
\buildrel{\text{law}}\over=
\Big
(L(T_U(k+1,x),y), \ \ y\in\Z
\Big).
\tag 3.6
$$



\bigskip
\noindent
{\bf 4. Proof of Theorem 1.}
 
\medskip
\noindent
Given the Markov chains $Y_t$, $Z_t$ and $Y'_t$ we define
$$
\wt Z_t
:=
Z_t + Z_{t-1} + 1, 
\qquad
\wt Y_t
:=
Y_t + Y_{t-1},
\qquad
\wt Y'_t
:=
Y'_t + Y'_{t-1}
\tag 4.1
$$ 
and for $h\in\N$ the following stopping times
$$
\align
\sigma_{h} 
&
:=
\inf\{t\ge0:Y_t\ge h\},
\tag 4.2
\\
\sigma'_{h} 
&
:=
\inf\{t\ge0:Y'_t\ge h\},
\tag 4.3
\\
\omega 
&
:=
\inf\{ t\ge 0 : Y_t =  0 \},
\tag 4.4
\\
\omega' 
&
:=
\inf\{ t\ge 0 : Y'_t =  0 \},
\tag 4.5
\\
\tau_{h} 
&
:=
\inf\{t\ge0:Z_t\ge h\},
\tag 4.6
\\
\wt\sigma_{h,0} 
&
:=
0,
\quad
\wt\sigma_{h,i+1}
:=
\inf\{t>\wt\sigma_{h,i}:\wt Y_t\ge h\},
\quad
\wt\sigma_h
:= 
\wt\sigma_{h,1},
\tag 4.7
\\
\wt\sigma'_{h,0} 
&
:=
0, 
\quad
\wt\sigma'_{h,i+1}
:=
\inf\{t>\wt\sigma'_{h,i}:\wt Y'_t\ge h\},
\quad
\wt\sigma'_h
:= 
\wt\sigma'_{h,1},
\tag 4.8
\\
\wt\tau_{h,0} 
&
:=
0, 
\quad
\wt\tau_{h,i+1}
:=
\inf\{t>\wt\tau_{h,i}:\wt Z_t\ge h\},
\quad
\wt\tau_h
:= 
\wt\tau_{h,1}.
\tag 4.9
\endalign
$$
For $h\ge1$, $p\ge0$ and $x\ge1$ fixed integers 
we define the 
following events:
$$
\align
A_{h,0}
:=
&
\left\{ 
\max \{\wt Y_t:  1 \le t < \infty\} < h 
\right\}
\\
=
&
\left\{ 
\phantom{\wt Y_t}
\wt\sigma_h = \infty 
\phantom{\wt Y_t}
\right\},
\tag 4.10 
\\
A_{h,p}
:=
&
\left\{ 
[\max \{\wt Y_t: 1\le t < \infty\}=h]
\ \land  \ 
[\#\{1 \le t < \infty: \wt Y_t=h\} =p] 
\right\}
\\
=
&
\left\{ 
[\wt\sigma_{h,p}<\infty=\wt\sigma_{h,p+1}]
\ \land \ 
[\wt Y_{\wt\sigma_{h,i}}=h, \ i=1,\dots,p] 
\right\},
\tag 4.11
\\
A'_{h,0}
:=
&
\left\{ 
\max\{ \wt Y'_t: 1\le t < \infty\} < h 
\right\}
\\
=
&
\left\{ 
\phantom{\wt Y_t}
\wt\sigma'_h = \infty 
\phantom{\wt Y_t}
\right\},
\tag 4.12
\\
A'_{h,p}
:=
&
\left\{ 
[\max \{\wt Y'_t: 1\le t < \infty\}=h] 
\ \land  \ 
[\#\{1 \le t < \infty: \wt Y'_t=h\} =p]
\right\}
\\
=
&
\left\{ 
[\wt\sigma'_{h,p}<\infty=\wt\sigma'_{h,p+1}]
\ \land \ 
[\wt Y'_{\wt\sigma'_{h,i}}=h, \ i=1,\dots,p] 
\right\},
\tag 4.13
\endalign
$$
$$
\align
B_{x,h,0}
:=
&
\left\{ 
\max \{\wt Z_t:1\le t < x\} < h 
\right\}
\\
=
&
\left\{ 
\phantom{\wt Y_t}
\wt\tau_h \ge x 
\phantom{\wt Y_t}
\right\},
\tag 4.14
\\
B_{x,h,p}
:=
&
\left\{ 
[\max \{\wt Z_t:1\le t <x\}=h] 
\ \land  \ 
[\#\{1\le t<x: \wt Z_t=h\} =p]
\right\}
\\
=
&
\left\{ 
[\wt\tau_{h,p}<x\le\wt\tau_{h,p+1}] 
\ \land \ 
[\wt Z_{\wt\tau_{h,i}}=h, \ i=1,\dots,p] 
\right\}.
\tag 4.15
\endalign
$$
With the help of the Ray-Knight representation and the events
introduced above we get the expression:
$$
\align
\expect\Big(u(x)\Big) = 
\sum_{p+q+r=3}
\sum_{h=1}^\infty
\sum_{k=0}^\infty
\sum_{l=0}^\infty
\,\,\,\,
&
\prob\Big(
A_{h,p}
\Big| Y_0=h-k-1 \Big) 
\times
\\
&
\pi(k,h-k-1) 
\times
\\
&
\prob\Big(
B_{x,h,q} \,\land\, [Z_{x-1}=l]
\Big| Z_0=k \Big)
\times
\\
&
\prob\Big(
A'_{h,r}
\Big| Y'_0=l \Big)
\tag 4.16
\endalign
$$
which leads directly to
$$
\align
\sum_{x=1}^\infty \expect\Big(u(x)\Big) \le
\sum_{p+q+r=3}
\sum_{h=1}^\infty
\sum_{k=0}^\infty
&
\phantom{M}
\prob\Big(A_{h,p}\Big|Y_0=h-k-1\Big)
\times
\\
&
\phantom{M}
\pi(k,h-k-1)
\times
\\
&
\left(
\sum_{x=1}^\infty
\prob\Big(B_{x,h,q}\Big|Z_0=k\Big)
\right)
\times
\\
&
\left(
\sup_{l\ge0}
\prob\Big(A'_{h,r}\Big|Y'_0=l\Big)
\right). 
\tag 4.17
\endalign
$$

The proof of Theorem 1 will follow directly from the bounds
provided by

\proclaim{Proposition 1}
For any $\vareps>0$ there exists a finite constant $C<\infty$
such that for any $h\ge1$ and $k\ge0$:
\item{($\imath$)} 
Without any restriction on $k$ or $p$
$$
\sum_{x=1}^\infty
\prob\Big( B_{x,h,p} \Big| Z_0=k \Big) 
\le
C h
\tag 4.18
$$
\item{($\imath\imath$)}
if either $k\in[(h- h^{1/2+\vareps})/2,(h- h^{1/2+\vareps})/2]$ 
or $p\ge 1$ holds then
$$
\align
\prob\Big( A_{h,p} \Big| Y_0=k \Big) 
\le
&
\left(C h^{-1/2+\vareps}\right)^{p+1}
\tag 4.19
\\
\sum_{x=1}^\infty
\prob\Big( B_{x,h,p} \Big| Z_0=k \Big) 
\le
&
\left(C h^{-1/2+\vareps}\right)^{p+1} h
\tag 4.20
\endalign
$$
\endproclaim

\noindent
{\sl Remark:} 
For $k\ge h$ the left hand side of (4.18), (4.19) and (4.20), 
of course, vanish.

We postpone the proof of this Proposition to the end  of the next 
section and proceed with the proof of Theorem 1.
Using the bounds (4.18)-(4.20) we prove (2.10). As
we already mentioned, (2.11) is proved in a completely identical
way. Theorem 1 follows from (2.10), (2.11) via (2.7). 

In the forthcoming argument $\cdots$ will stand as 
abreviation of the summand on the right hand side of (4.17).
On the right hand side of (4.17) keep $p,q,r$ and $h$ fixed and
decompose the sum over $k\ge0$ as follows:
$$
\sum_{k}
\cdots
\,\,\,\,
=
\sum_{k:|h-2k|\le h^{1/2+\vareps}}
\!\!\!\!\!\!\!\!\!\!\!\!
\cdots
\,\,\,\,
+
\sum_{k:|h-2k|> h^{1/2+\vareps}}
\!\!\!\!\!\!\!\!\!\!\!\!
\cdots
\,.
\tag 4.21
$$
Similar decompositions will be applied a few more times 
throughout the paper. 

\proclaim{Side-lemma 1}
For any $\vareps>0$ there exist  constants $C<\infty$ and
$\gamma>0$  such that for any $h\ge 1$
$$
\sum_{k:|h-2k|> h^{1/2+\vareps}}
\!\!\!\!\!\!\!\!\!\!\!\!
\pi(k, h-1-k)
<
C \exp(-\gamma h^{2\vareps}). 
\tag 4.22
$$
\endproclaim

\noindent 
Side-lemmas are proved in Section 6.

Using  (4.18)-(4.21) we bound the sum over $k$ on the right hand
side of (4.17) as follows: 
\newline
If $r=0$ and  $p+q=3$
$$
\align
\sum_{k}
\cdots
\le 
&
\,\,\,
\left( C h^{-1/2+\vareps} \right)^{p+q+2} h
+
\big( C h \big)
\big( C \exp(-\gamma h^{2\vareps}) \big)
\\
\le 
&
\,\,\,
C' h^{-3/2+5\vareps}
\tag 4.23
\endalign
$$
with some properly chosen $C'<\infty$. 
\newline
If $r>0$ and  $p+q+r=3$
$$
\align
\sum_{k}
\cdots
\le 
&
\,\,\,
\left(C h^{-1/2+\vareps}\right)^{p+q+r+3} h
+
\big( C h \big)
\big( C \exp(-\gamma h^{2\vareps}) \big)
\\
\le
&
\,\,\,
C' h^{-2+6\vareps}
\tag 4.24
\endalign
$$
with some properly chosen $C'<\infty$. 

In both cases the upper bound is summable over $h\ge1$, if we
choose $\vareps<1/10$. Hence (2.10) and the statement of Theorem
1.  

\noindent
\qed (Theorem 1)


\bigskip
\noindent
{\bf 5. Technical Lemmas.}

\medskip 
\noindent
The present section is divided in five subsections. In subsections 
5.1-5.4 we state and prove some lemmas of more technical nature, needed 
in the proof of Proposition 1, which is presented in subsection 5.5.
Throughout this section $\vareps>0$ is fixed.
 
\medskip
\noindent
{\sl 5.1. The Maximal Jump.}

\smallskip
\noindent
We prove that the largest jump of the Markov chains $Y_t$ and
$Z_t$, before reaching level $h$, is less than
$h^{1/2+\vareps}$, with overwhelming probability. 
Define the maximal jumps of $Y_t$, respectively, $Z_t$ as
follows: 
$$
\align
M_h
:=
&
\sup
\left\{
|Y_t-Y_{t-1}|: 1\le t \le \sigma_h
\right\}
\\
=
&
\sup
\left\{
|Y_t-Y_{t-1}|: 1\le t \le \sigma_h \land \omega
\right\},
\tag 5.1 
\\
N_h
:=
&
\sup
\left\{
|Z_t-Z_{t-1}|: 1\le t \le \tau_h
\right\}.
\tag 5.2
\endalign
$$
By definition $M_h=0$ if $Y_0\ge h$ and $N_h=0$ if $Z_0\ge h$.  

\proclaim{Lemma 1}
There exist two constants, $C<\infty$ and $\gamma>0$, such that
for any $h\ge1$ and $k\ge0$ the following bounds hold:
$$
\align
\prob\Big(M_h>h^{1/2+\vareps}\Big|Y_0=k\Big)
< C\exp(-\gamma h^{2\vareps}),
\tag 5.3
\\
\prob\Big(N_h>h^{1/2+\vareps}\Big|Z_0=k\Big)
< C\exp(-\gamma h^{2\vareps}).
\tag 5.4
\endalign
$$
\endproclaim

\demo{Proof}
We prove here (5.3) in details. The proof of (5.4) is
essentially similar and it is left for the reader. For the
moment let $\gamma$ be an arbitrary positive number. Its value will
be fixed at the end of this proof. 
$$
\align
\prob\Big( M_h>h^{1/2+\vareps} 
\,\,
&
\Big| Y_0=k \Big) 
\\
\le
\phantom{+}
&
\prob\Big( 
[M_h>h^{1/2+\vareps}] 
\ \land \  
[\sigma_h\,\land\,\omega \le h^2\exp(\gamma h^{2\vareps})]
\Big| Y_0=k \Big)
\\
+
&
\prob\Big(
\sigma_h\,\land\,\omega > h^2\exp(\gamma h^{2\vareps})
\Big| Y_0=k \Big)
\tag 5.5
\endalign
$$

\proclaim{Side-lemma 2}
There exists a constant $C<\infty$, such that for any $h\ge1$
and $k\ge0$ 
$$
\expect\Big( 
\sigma_h\,\land\,\omega
\Big| Y_0=k\Big) 
\le C h^2.
\tag 5.6
$$
\endproclaim

\noindent
By using Markov's inequality we get the following upper bound on
the second term of the right hand side of (5.5):
$$
\prob\Big(
\sigma_h \,\land\, \omega >
h^2 \exp( \gamma h^{2\vareps} ) 
\Big| Y_0=k \Big)
\le  
C \exp( -\gamma h^{2\vareps} ).
\tag 5.7
$$
To bound the first term on the right hand side of (5.5) we use
the following representation of the Markov chain $Y_t$:
Let $\left( \xi_{t,i} \right)_{t\ge 1, i\ge 1}$ 
be i.i.d random variables
with common geometric distribution 
$\prob\Big(\xi_{t,i}=k\Big)=2^{-k-1}$
The process $Y_t$ is realized as follows: fix $Y_0$ and put
$$
Y_{t+1}=\sum_{j=1}^{Y_t}\xi_{t+1,j}.
\tag 5.8
$$
Using this representation we note that
$$
\align
&
\prob\Big( 
[M_h > h^{ 1/2 + \vareps }] 
\ \land \  
[\sigma_h \,\land\, \omega \le h^2\exp(\gamma h^{2\vareps})]
\Big| Y_0=k \Big)
\\
&
\qquad\qquad
\le 
\prob\Big(
\max
\left\{
\max_{ 1\le j \le h } 
\left|\sum_{i=1}^{j} 
\left( \xi_{t,i} - 1 \right)
\right|:
1 \le t \le h^2 \exp(\gamma h^{2\vareps})
\right\}
> h^{ 1/2 + \vareps } 
\Big)
\\
&
\qquad\qquad
=
1-
\left(
1-
\prob\Big(
\max_{ 1\le j \le h } 
\left|
\sum_{i=1}^{j}
\left( \xi_{1,i} - 1 \right)
\right|
> h^{ 1/2 + \vareps }
\Big) 
\right)^{ h^2 \exp(\gamma h^{2\vareps}) }
\tag 5.9
\endalign
$$

\proclaim{Side-lemma 3}
Let $\xi_i$ be i.i.d. random variables with the common geometric
disatribution $\prob\Big(\xi_i=k\Big)=2^{-k-1}$. Then there is a 
constant $\theta_0>0$ such that for any
$\lambda>0$ and $n\in\N$ satisfying $\lambda/(4n)<\theta_0$:
$$
\prob\Big(
\max_{1\le j \le n} 
\left|
\sum_{i=1}^{j}
\left( \xi_{i}-1\right)
\right|
> \lambda
\Big) 
\le 2 \exp(-\lambda^2/(8n)). 
\tag 5.10
$$
\endproclaim

\noindent
Using this bound we get
$$
\align
&
\prob\Big(
[M_h > h^{ 1/2 + \vareps }]  
\ \land \ 
[\sigma_h \,\land\, \omega \le h^2 \exp(\gamma h^{2\vareps})]
\Big| Y_0=k \Big)
\\
&
\qquad\qquad\qquad\qquad\qquad\qquad\qquad
\le
1-
\left(1-2\exp(-h^{2\vareps}/8)\right)^{h^2\exp(\gamma h^{2\vareps})}
\\
&
\qquad\qquad\qquad\qquad\qquad\qquad\qquad
\le
2h^2\exp\Big((\gamma-8^{-1})h^{2\vareps}\Big).
\tag 5.11
\endalign
$$
In the last inequality we use the fact 
that for $0<\alpha<1 <\beta$, \ 
$1-\alpha\beta <(1-\alpha)^\beta$.
We choose  $\gamma<16^{-1}$. From (5.5), (5.7) and (5.11) we get 
(5.3), with an
appropriately chosen constant $C<\infty$.

\noindent
\qed (Lemma 1)

\medskip
\noindent
{\sl 5.2. Hitting \underbar{exactly} $h$.}

\smallskip
\noindent

\proclaim{Lemma 2}
There exists a constant $C<\infty$ such that for any $h\ge 1$
and $k\ge 0$ 
$$
\align
\prob\Big(
[\wt\sigma_h<\infty] 
\ \land \ 
[\wt Y_{\wt\sigma_h}=h]
\Big| Y_0=k \Big) 
< 
C h^{-1/2+\vareps}
\tag 5.12
\\
\prob\Big( 
\wt Z_{\wt\tau_h}=h
\Big| Z_0=k \Big) 
< 
C h^{-1/2+\vareps}
\tag 5.13
\endalign
$$
\endproclaim

\demo{Proof}
Again, we give the details of the proof of (5.12), leaving the 
identical details  of (5.13) for the reader.  
$$
\align
&
\prob\Big(
[\wt\sigma_h<\infty] 
\ \land \  
[\wt Y_{\wt\sigma_h}=h]
\Big| Y_0=k \Big) 
=
\\
&
\qquad\qquad\qquad\qquad
\sum_{l=0}^\infty
\prob\Big(
[\wt\sigma_h<\infty] 
\ \land \ 
[Y_{\wt\sigma_h-1}=l] 
\ \land \ 
[Y_{\wt\sigma_h}=h-l]  
\Big| Y_0=k \Big)
\tag 5.14
\endalign
$$
We divide the sum in two parts, as in (4.21):
$$
\align
\sum_{l:|h-2l|>h^{1/2+\vareps}}
\!\!\!\!\!\!\!\!\!
&
\prob\Big(
[\wt\sigma_h<\infty]  
\ \land \ 
[Y_{\wt\sigma_h-1}=l]
\ \land \ 
[Y_{\wt\sigma_h}=h-l]  
\Big| Y_0=k \Big) 
\\
&
\qquad\qquad
\le
\prob\Big(
M_h>h^{1/2+\vareps} 
\Big| Y_0=k \Big)
<
C\exp(-\gamma h^{2\vareps}),
\tag 5.15
\endalign
$$
by Lemma 1. On the other hand:
$$
\align
&
\sum_{l:|h-2l|\le h^{1/2+\vareps}}
\!\!\!\!\!\!\!\!\!
\prob\Big(
[\wt\sigma_h<\infty] 
\ \land \ 
[Y_{\wt\sigma_h-1}=l] 
\ \land \ 
[Y_{\wt\sigma_h}=h-l]  
\Big| Y_0=k \Big) 
=
\\
&
\quad\quad\quad
\sum_{l:|h-2l|\le h^{1/2+\vareps}}
\!\!\!\!\!\!\!\!\!
\prob\Big(
[\wt\sigma_h<\infty]
\ \land \ 
[Y_{\wt\sigma_h-1}=l]
\Big| Y_0=k \Big) 
\frac{\pi(l,h-l)}{\sum_{m\ge h-l} \pi(l,m)}
\tag 5.16
\endalign
$$

\proclaim{Side-lemma 4}
There exists a constant $C<\infty$, such that for any
$h\ge1$   and $l\in[ (h-h^{1/2+\vareps})/2,  (h+h^{1/2+\vareps})/2]$
$$
\frac{\pi(l,h-l)}{\sum_{m\ge h-l} \pi(l,m)}
<
C h^{-1/2+\vareps}
\tag 5.17
$$
\endproclaim

\noindent
From this we get 
$$
\sum_{l:|h-2l|\le h^{1/2+\vareps}}
\!\!\!\!\!\!\!\!\!
\prob\Big(
[\wt\sigma_h<\infty] 
\ \land \ 
[Y_{\wt\sigma_h-1}=l] 
\ \land \ 
[Y_{\wt\sigma_h}=h-l]  
\Big| Y_0=k \Big) 
\le 
C h^{-1/2+\vareps}.
\tag 5.18
$$
Finally, (5.15) and (5.18) yield (5.13).

\noindent
\qed (Lemma 2)
\enddemo

\medskip
\noindent
{\sl 5.3. $\wt Y_t$ does not hit level $\ge h$.}

\smallskip
\noindent

\proclaim{Lemma 3}
There exists a constant $C<\infty$ such that for any $h\ge1$ 
and $k\in [ (h-h^{1/2+\vareps})/2 , (h+h^{1/2+\vareps})/2]$
$$
\prob\Big(
\wt\sigma_h=\infty
\Big| Y_0=k \Big)
<
C h^{-1/2+\vareps}.
\tag 5.19
$$
\endproclaim

\demo{Proof}
$$
\align
\prob\Big(
\wt\sigma_h=\infty
\Big| Y_0=k \Big)
\le
\phantom{+}
&
\,\,
\prob\Big(
[\wt\sigma_h=\infty]
\ \land \ 
[M_h\le h^{1/2+\vareps}]
\Big| Y_0=k \Big) 
\\
+
&
\,\,
\prob\Big(
M_h > h^{1/2+\vareps} 
\Big| Y_0=k \Big).
\tag 5.20
\endalign
$$
To bound the first term on the right hand side, note that
$$
\left\{ 
[\wt\sigma_h=\infty]
\ \land \ 
[M_h\le h^{1/2+\vareps}]
\right\}
\subset
\left\{
\phantom{h^h_h}
\sigma_{(h+h^{1/2+\vareps})/2}=\infty
\phantom{h^h_h}
\right\}.
\tag 5.21
$$

\proclaim{Side-lemma 5}
There exists a constant $C<\infty$ such that for any $0\le k < h$
$$
\prob\Big(
\sigma_{h}=\infty
\Big| Y_0=k \Big)
<
\frac{h-k}{h}+ C h^{-1/2}.
\tag 5.22
$$
\endproclaim

\noindent
Thus
$$
\align
\prob\Big(
[\wt\sigma_h=\infty]
\ \land \ 
[M_h \le h^{1/2+\vareps}]
\Big| Y_0=k \Big)
\le
&\,\,
\prob\Big(
\sigma_{(h+h^{1/2+\vareps})/2}=\infty
\Big| Y_0=k \Big)
\\
\le
&\,\,
C h^{-1/2 +\vareps},
\tag 5.23
\endalign
$$
for $k\in [ (h-h^{1/2+\vareps})/2 , (h+h^{1/2+\vareps})/2]$.
This bound, together with (5.20) and (5.3) yield (5.19).

\noindent
\qed (Lemma 3)
\enddemo

\medskip
\noindent
{\sl 5.4. Expectation of $\wt\tau_h$.}

\smallskip
\noindent

\proclaim{Lemma 4}
There exists a  constant  $C <\infty$
such that for any $h\ge1$  the following bounds hold:
\item{($\imath$)}
Without any restriction on $k$:
$$
\expect\Big(
\wt\tau_h
\Big| Z_0=k \Big)
<
C h.
\tag 5.24
$$
\item{($\imath\imath$)}
For $k\in[(h-h^{1/2+\vareps})/2, (h+h^{1/2+\vareps})/2]$
$$
\expect\Big(
\wt\tau_h
\Big| Z_0=k \Big)
<
C h^{1/2+\vareps}.
\tag 5.25
$$
\endproclaim
 
\demo{Proof}
$$
\expect\Big(
\wt\tau_h
\Big| Z_0=k \Big)
=
\expect\Big(
\wt\tau_h \ind_{\{N_h\le h^{1/2+\vareps}\}}
\Big| Z_0=k \Big)
+
\expect\Big(
\wt\tau_h \ind_{\{N_h> h^{1/2+\vareps}\}}
\Big| Z_0=k \Big).
\tag 5.26
$$
We bound the first, respectively, the second term on the right
hand side, by noting 
$$
\wt\tau_h \ind_{\{N_h\le h^{1/2+\vareps}\}} 
\le 
\tau_{(h+h^{1/2+\vareps})/2},
\tag 5.27
$$
respectively,
$$
\wt\tau_h
\le
\tau_h.
\tag  5.28
$$
Thus we get
$$
\align
\expect\Big(
\wt\tau_h 
\Big| Z_0=k \Big)
\le
\phantom{+}
&
\phantom{+}
\expect\Big(
\tau_{(h+h^{1/2+\vareps})/2} 
\Big| Z_0=k \Big)
\\
&
+
\expect\Big(
\tau_h^2 
\Big| Z_0=k \Big)^{1/2}
\prob\Big(
N_h> h^{1/2+\vareps} 
\Big| Z_0=k \Big)^{1/2}.
\tag 5.29
\endalign
$$

\proclaim{Side-lemma 6}
There exists a constant $C<\infty$ such that for any 
$0\le k < h$
$$
\expect\Big(
\tau_{h} 
\Big| Z_0=k \Big)
<
(h-k) + C h^{1/2}.
\tag 5.30
$$
\endproclaim

\proclaim{Side-lemma 7}
There exists a constant $C <\infty$ such that for any 
$0\le k < h$
$$
\expect\Big(
\tau_{h}^2 
\Big| Z_0=k \Big)
<
C h^2.
\tag 5.31
$$
\endproclaim

\noindent
Putting together (5.29), (5.30) and (5.31), we get 
(5.24) and (5.25).

\noindent
\qed (Lemma 4)
\enddemo

\medskip
\noindent
{\sl 5.5 Proof of Proposition 1.}

\smallskip
\noindent
First note that
$$
\align
\prob\Big(
A_{h,0}
\Big| Y_0=k \Big)
& 
=
\prob\Big(
\wt\sigma_h=\infty
\Big| Y_0=k \Big),
\tag 5.32
\\
\sum_{x=1}^\infty
\prob\Big(
B_{x,h,0}
\Big| Z_0=k \Big)
& 
=
\sum_{x=1}^\infty
\prob\Big(
\wt\tau_h\ge x
\Big| Z_0=k \Big)
=
\expect\Big(
\wt\tau_h
\Big| Z_0=k \Big),
\tag 5.33
\endalign
$$
and for $p\ge1$, using the strong Markov property of $Y_t$,
respectively, $Z_t$:
$$
\align
\prob\Big(
&
A_{h,p}
\Big| Y_0=k \Big)
=
\sum_{l=0}^\infty
\prob\Big(
[\wt\sigma_h<\infty]
\ \land \ 
[Y_{\wt\sigma_h-1}=h-l]
\ \land \
[Y_{\wt\sigma_h}=l]
\Big| Y_0=k \Big)
\times
\\
&
\qquad\qquad\qquad\qquad\quad
\prob\Big(
A_{h,p-1}
\Big| Y_0=l \Big)
\tag 5.34
\\
\sum_{x=1}^\infty
&
\prob\Big(
B_{x,h,p}
\Big| Z_0=k \Big)
=
\sum_{l=0}^\infty
\prob\Big(
[Z_{\wt\tau_h-1}=h-l]
\ \land \
[Z_{\wt\sigma_h}=l]
\Big| Z_0=k \Big)
\times
\\
&
\qquad\qquad\qquad\qquad\qquad\quad
\left(
\sum_{x=1}^\infty
\prob\Big(
B_{x,h,p-1}
\Big| Z_0=l \Big)
\right).
\tag 5.35
\endalign
$$
We prove the bounds of Proposition 1 by induction on $p$. 

According to (5.32), (5.33), 
for $p=0$, 
(4.18), (4.19) and (4.20) 
are just restatements of 
(5.24), (5.19) and (5.25), 
respectively. (See Lemma 3 and Lemma 4.)

Next we consider the case $p=1$. 
Again, we divide the sum over $l$ in (5.34) and (5.35) 
in two parts, as it was done in (4.21). 
From (5.12) (Lemma 2) and (5.19) (Lemma 3)
$$
\align
\sum_{l:|h-2l|\le h^{1/2+\vareps}}
&
\prob\Big(
[\wt\sigma_h<\infty]
\ \land \ 
[Y_{\wt\sigma_h-1}=h-l]
\ \land \
[Y_{\wt\sigma_h}=l]
\Big| Y_0=k \Big)
\prob\Big(
A_{h,0}
\Big| Y_0=l \Big)
\\
&
\qquad\qquad\qquad
\le 
\left( C  h^{-1/2+\vareps}\right)
\left( C h^{-1/2+\vareps}\right).
\tag 5.36
\endalign
$$
From (5.3) (Lemma 1)
$$
\align
&
\sum_{l:|h-2l| > h^{1/2+\vareps}}
\!\!\!\!\!\!\!\!\!
\prob\Big(
[\wt\sigma_h<\infty]
\ \land \ 
[Y_{\wt\sigma_h-1}=h-l]
\ \land \
[Y_{\wt\sigma_h}=l]
\Big| Y_0=k \Big)
\prob\Big(
A_{h,0}
\Big| Y_0=l \Big)
\\
&
\qquad\qquad\quad\qquad\qquad\quad
\le 
\prob\Big(
M_h>h^{1/2+\vareps}
\Big| Y_0=k \Big)
<
C\exp(-\gamma h^{2\vareps})
\tag 5.37
\endalign
$$
From (5.36) and (5.37) we get (4.19) for $p=1$. 

Applying the same ideas to (5.35): from (5.13) (Lemma 2) 
and (5.25) (Lemma 4) 
$$
\align
&
\sum_{l: |h-2l|\le h^{1/2+\vareps}}
\!\!\!\!\!\!\!\!\!
\prob\Big(
[Z_{\wt\tau_h-1}=h-l]
\ \land \
[Z_{\wt\sigma_h}=l]
\Big| Z_0=k \Big)
\sum_{x=1}^\infty
\prob\Big(
B_{x,h,0}
\Big| Z_0=l \Big)
\\
&
\qquad\qquad\qquad\qquad\qquad
\le 
\left(C h^{-1/2+\vareps}\right)
\left(C h^{1/2+\vareps}\right)
\tag 5.38
\endalign
$$
From (5.4) (Lemma 1) and (5.24) (Lemma 4)
$$
\align
&
\sum_{l:|h-2l| > h^{1/2+\vareps}}
\!\!\!\!\!\!\!\!\!
\prob\Big(
[Z_{\wt\sigma_h-1}=h-l]
\ \land \
[Z_{\wt\sigma_h}=l]
\Big| Z_0=k \Big)
\sum_{x=1}^\infty
\prob\Big(
B_{x,h,0}
\Big| Z_0=l \Big)
\\
&
\qquad\qquad\qquad
\le
\prob\Big(
N_h>h^{1/2+\vareps}
\Big| Z_0=k \Big)
\left(
\sup_{l\ge0}
\sum_{x=1}^\infty
\prob\Big(
B_{x,h,0}
\Big| Z_0=l \Big)
\right)
\\
&
\qquad\qquad\qquad
<
\left( C\exp(-\gamma h^{2\vareps}) \right)
(C h )
\tag 5.39
\endalign
$$
(5.38) and (5.39) yield (4.20) for $p=1$. 

For $p\ge2$ the induction follows from the same reasonings, 
just one does not have to split the sum over $l$ as in (5.36), 
(5.37). After the previous arguments we may ignore 
these completely straightforward details.
 
\noindent
\qed (Proposition 1)



\bigskip
\noindent
{\bf 6. Proof of the Side-Lemmas.}
 
\medskip
\noindent
First we prove Side-lemmas 1 and 4.  
Then Side-lemma 3 follows, which 
relies on an exponential Kolmogorov 
inequality. These proofs are rather 
standard `classroom exercises'. 
Side-lemmas 2, 5, 
6 and 7 follow from an estimate on the 
overshooting of level $h$ 
by the processes $Y_t$ and $Z_t$ stopped 
at $\sigma_h \land \omega$, 
respectively, $\tau_h$ and from optional 
stopping arguments. 

\demo{Proof of Side-lemma 1}
\newline
Assume $h\ge2$ and denote $h-2=:n$, $k-1=:l$. Then, using the
explicit form (3.1) of $\pi(i,j)$, the right hand side of (4.22)
becomes 
$$
\sum_{k:|h-2k|>h^{1/2+\vareps}}
\!\!\!\!\!\!\!\!\!
\pi(k,h-1-k)
=
\frac12
\prob\Big(
\left|2B_n - n\right| > (n+2)^{1/2+\vareps}
\Big)
\tag 6.1
$$
where $B_n$ is binomially distributed: 
$\prob\left( B_n=l \right)={n\choose l}2^{-n}$. Using the fact
that for any $\gamma<1/2$ 
$$
\sup_n\expect\Big(
\exp\left\{ \gamma \left( 2B_n-n \right)^2 / n \right\} 
\Big)
= 
C_\gamma<\infty,
\tag 6.2
$$
by Markov's inequality  we get (4.22).

\noindent
\qed
(Side-lemma 1)
\enddemo

\demo{Proof of Side-lemma 4}
\newline
Note first that for $i\ge1$ and $j\ge0$
$$
\frac{\pi(i,j+1)}{\pi(i,j)}
=
\frac{i+j}{2(1+j)}
\tag 6.3
$$
From this it follows that the distribution
$j\mapsto\pi(i,j)$ is unimodular and 
for $i\ge2$ fixed 
$$
\aligned
&
\pi(i,j)
<\pi(i,j+1),
\quad\quad\quad
\text{ for } 0\le j \le i-3,
\\
&
\pi(i,i-2) 
= 
\pi(i,i-1),
\\
&
\pi(i,j)
>
\pi(i,j+1)
\phantom{,}
\quad\quad\quad
\text{ for } i-1\le j <\infty. 
\endaligned
\tag 6.4
$$
We treat separately the cases 
$l\in[h/2,(h+h^{1/2+\vareps})/2]$ and 
$l\in[(h-h^{1/2+\vareps})/2, h/2]$:
\newline
For $l\in[h/2,(h+h^{1/2+\vareps})/2]$ the following two
facts imply (5.17)
\newline
(1) 
By (6.4),  
$$
\align
\pi(l,h-l)
&
\le 
\pi(l,l-1)
=
\frac12 {{2(l-1)}\choose{l-1}}2^{-2(l-1)}
\\
&
\le
\frac12 
{{2(\lfloor h/2\rfloor-1)}\choose{\lfloor h/2\rfloor-1}} 
2^{-2(\lfloor h/2\rfloor-1)}
\le
C h^{-1/2}
\tag 6.5
\endalign
$$
(2) By the central limit theorem: 
$\lim_{l\to\infty}\sum_{m\ge l} \pi(l,m) =\frac12$ and thus
there exists a constant $c>0$ such that for any $0<h/2\le
l$ 
$$ 
\sum_{m\ge h-l} \pi(l,m) \ge \sum_{m\ge l} \pi(l,m) \ge c
\tag 6.6
$$
Let now 
$l\in[(h+h^{1/2-\vareps})/2, h/2]$ and 
$k := \lfloor h-l + h^{1/2-\vareps}\rfloor$
Then:
$$
\align
\frac{\pi(l,h-l)}{\sum_{m\ge h-l}\pi(l,m)}
&
\le
(k-h+l+1)^{-1} \frac{\pi(l,h-l)}{\pi(l,k)}
\\
&
\le
(k-h+l+1)^{-1} 
\left(\frac{\pi(l,k-1)}{\pi(l,k)}\right)^{k-h+l}
= 
\cdots
\endalign
$$
$$
\align
\phantom{\frac{\pi(l,h-l)}{\sum_{m\ge h-l}\pi(l,m)}}
\cdots
&
=
(k-h+l+1)^{-1} 
\left(\frac{2k}{l+k-1}\right)^{k-h+l}
\\
&
\le
h^{-1/2+\vareps}
\left(
\frac{2(h-l + h^{1/2-\vareps})} 
     {h+ h^{1/2-\vareps}-1}
\right)^{h^{1/2-\vareps}}
\\
&
\le
h^{-1/2+\vareps}
\left(
\frac{h+ h^{1/2+\vareps}+ 2h^{1/2-\vareps})} 
     {h+ h^{1/2-\vareps}-1}
\right)^{h^{1/2-\vareps}}
\\
&
\le
h^{-1/2+\vareps}
\left( 1 + 3 h^{-(1/2-\vareps)}\right)^{h^{1/2-\vareps}}
\\
&
\le 
e^3 h^{-1/2+\vareps}
\tag 6.7
\endalign
$$
In the first inequality we used (6.1). In the second one we 
exploited the fact that, according to (6.2), for any $i\ge1$
fixed $\pi(i,j)/\pi(i,j+1)$ is an increasing function of $j\ge0$. 
In the next equality, (6.2) was explicitly used. In the third 
inequality we inserted the value of 
$k=\lfloor h-l+h^{1/2-\vareps}\rfloor$. In 
the fourth inequality 
$l=\lceil (h-h^{1/2-\vareps})/2\rceil$ was inserted to 
maximize the expression. In  the last but one inequality we used 
$1\le h^{1/2-\vareps}\le h^{1/2+\vareps}$. Finally, in the last 
inequality we use the fact that 
$\sup_{\alpha\ge1}\left(1+3\alpha^{-1}\right)^\alpha\le e^3$

\noindent
\qed
(Side-lemma 4)
\enddemo

\demo{Proof of Side-lemma 3}
\newline
The Exponential Kolmogorov Inequality follows
directly from Doob's maximal inequality. For its proof  
see e.g. page 139 of  Williams (1991).

\proclaim{Exponential Kolmogorov Inequality}
Let $\wt\xi_j$, $j\ge1$, be i.i.d. random variables with 
$\expect\Big(
\exp\left\{\theta \left| \wt\xi_j \right|\right\}
\Big)<\infty$ for some $\theta>0$ and  
$\expect\Big(\wt\xi_j\Big)=0$. Then for any $\lambda\in(0,\infty)$
and $n\in\N$
$$
\prob\Big(
\max_{1\le j\le n}\left| \sum_{i=1}^j \wt\xi_i \right| >\lambda
\Big)
\le
e^{-\lambda\theta}
\left\{
\expect\Big(e^{\theta \wt\xi_i}\Big)^n
+
\expect\Big(e^{-\theta \wt\xi_i}\Big)^n
\right\}
\tag 6.8
$$
\endproclaim

\noindent
We apply the Exponential Kolmogorov Inequality
to $\wt\xi_i=\xi_i-1$, with 
$\prob\Big(\xi_i=k\Big)=2^{-k-1}$, $k\ge0$. 
There exists a constant $\theta_0>0$ such that for 
$0\le\theta<\theta_0$ we get:
$$
\align
&
\expect\Big(
e^{\theta (\xi_j-1)}
\Big)
=
e^{-\theta}\left(2-e^\theta\right)^{-1}
=
1+\theta^2+\Cal O\left(\theta^3\right)
<
e^{2\theta^2}
\tag 6.9
\\
&
\expect\Big(
e^{-\theta (\xi_j-1)}
\Big)
=
e^{2\theta}\left(2e^\theta-1\right)^{-1}
=
1+\theta^2+\Cal O\left(\theta^3\right)
<
e^{2\theta^2}
\tag 6.10
\endalign
$$
Inserting these bounds into the right hand side of (6.8) and choosing 
$\theta = \lambda/(4n)$ we obtain (5.10).

\noindent
\qed (Side-lemma 3)
\enddemo

The proof of Side-lemmas 2, 5, 6 and 7 will follow from the 
forthcoming  Overshooting Lemma and standard optional stopping 
considerations. The Overshooting Lemma and its 
Corollary are extended 
restatements of Lemmas 3.2 and 3.4 from T\'oth and Werner (1997).

\proclaim{Overshooting Lemma}
For any $0\le k < h \le u$ the following 
overshoot bounds hold:
$$
\align
\prob\Big(
Y_{\sigma_h}\ge u
\Big| [Y_0=k] \,\land\, [\sigma_h<\infty] \Big)
\le
&\,\,
\prob\Big(
Y_{1}\ge u
\Big| [Y_0=h] \,\land\, [Y_1\ge h] \Big)
\\
=
&\,\,
\frac{\sum_{v=u}^\infty \pi(h,v)}{\sum_{w=h}^\infty \pi(h,w)},
\tag 6.11
\\
\prob\Big(
Z_{\tau_h}\ge u
\Big| Z_0=k \Big)
\le
&\,\,
\prob\Big(
Z_{1}\ge u
\Big| [Z_0=h] \,\land\, [Z_1\ge h] \Big)
\\
=
&\,\,
\frac{\sum_{v=u}^\infty \rho(h,v)}{\sum_{w=h}^\infty \rho(h,w)}.
\tag 6.12
\endalign
$$
\endproclaim

In particular it follows that

\proclaim{Corollary}
There exists a constant $C<\infty$ such that for any $0\le k < h$:
$$
\align
\expect\Big(
Y_{\sigma_h}
\Big| [Y_0=k] \,\land\, [\sigma_h<\infty] \Big)
\le 
&\,\,
\frac{\sum_{v=h}^\infty \pi(h,v)v}{\sum_{w=h}^\infty \pi(h,w)}
\le 
h+C h^{1/2}
\tag 6.13
\\
\expect\Big(
Y_{\sigma_h}^2
\Big| [Y_0=k] \,\land\, [\sigma_h<\infty] \Big)
\le 
&\,\,
\frac{\sum_{v=h}^\infty \pi(h,v)v^2}{\sum_{w=h}^\infty \pi(h,w)}
\le 
h^2+C h^{3/2}
\tag 6.14
\\
\expect\Big(
Z_{\tau_h}
\Big| Z_0=k \Big)
\le 
&\,\,
\frac{\sum_{v=h}^\infty \rho(h,v)v}{\sum_{w=h}^\infty \rho(h,w)}
\le 
h+C h^{1/2}
\tag 6.15
\\
\expect\Big(
Z_{\tau_h}^2
\Big| Z_0=k \Big)
\le 
&\,\,
\frac{\sum_{v=h}^\infty \rho(h,v)v^2}{\sum_{w=h}^\infty \rho(h,w)}
\le 
h^2+C h^{3/2}
\tag   6.16
\endalign
$$
\endproclaim
The rightmost bounds in (6.13)-(6.16)
follow from explicit computations.

\demo{Proof of the Overshooting Lemma}
\newline
Straightforward manipulations leed to the following identities for $1\le h\le v$:
$$
\align
&
\prob\Big(
Y_{\sigma_h} = v
\Big| [Y_0=k] \,\land\, [\sigma_h<\infty] \Big)
=
\\
&
\qquad\qquad\quad
\sum_{l=0}^{h-1}
\prob\Big(
Y_{\sigma_h-1} =l
\Big| [Y_0=k] \,\land\, [\sigma_h<\infty] \Big)
\frac{\pi(l,v)}{\sum_{w=h}^\infty \pi(l,w)},
\tag 6.17
\endalign
$$
$$
\prob\Big(
Z_{\tau_h} = v
\Big| Z_0=k \Big)
=
\sum_{l=0}^{h-1}
\prob\Big(
Z_{\tau_h-1} =l
\Big| Z_0=k \Big)
\frac{\rho(l,v)}{\sum_{w=h}^\infty \rho(l,w)}.
\tag 6.18
$$
Using the explicit form (3.1), respectively, (3.2) of the transition 
probabilities $\pi(i,j)$, respectively, $\rho(i,j)$, it is easy 
to check the following inequalities, which hold for any 
$0<l<h\le v$, respectively, $0\le l<h\le v$:
$$
\align
\frac{\pi(l+1,v)}{\pi(l,v)}
=
\frac{l+v}{2l}
& 
<
\frac{l+v+1}{2l}
=
\frac{\pi(l+1,v+1)}{\pi(l,v+1)},
\tag 6.19
\\
\frac{\rho(l+1,v)}{\rho(l,v)}
=
\frac{l+v+1}{2(l+1)}
&
<
\frac{l+v+2}{2(l+1)}
=
\frac{\rho(l+1,v+1)}{\rho(l,v+1)}.
\tag 6.20
\endalign
$$
It follows that  for any $0\le l < h \le v < w$
$$
\align
&
\pi(l+1,v)\pi(l,w) 
<
\pi(l,v)\pi(l+1,w),
\tag 6.21
\\
&
\rho(l+1,v)\rho(l,w) 
<
\rho(l,v)\rho(l+1,w).
\tag 6.22
\endalign
$$
Hence, for any $0\le l< h \le u$
$$
\align
&
\sum_{v=h}^{\infty}\pi(l+1,v)
\sum_{w=u}^{\infty}\pi(l,w)
<
\sum_{v=h}^{\infty}\pi(l,v)
\sum_{w=u}^{\infty}\pi(l+1,w),
\tag 6.23
\\
&
\sum_{v=h}^{\infty}\rho(l+1,v)
\sum_{w=u}^{\infty}\rho(l,w)
<
\sum_{v=h}^{\infty}\rho(l,v)
\sum_{w=u}^{\infty}\rho(l+1,w),
\tag 6.24
\endalign
$$
which directly imply  (6.11), respectively, (6.12).

\noindent 
\qed (Overshooting Lemma)
\enddemo

\demo{Proof of Side-lemmas 2 and 5}
\newline
We apply the Optional Stopping Theorem to the martingales
$Y_t$ (for Side-lemma 5), 
respectively,
$Y_t^2-2\sum_{s=0}^{t-1}Y_s$ (for Side-lemma 2), 
$t\ge0$, 
both stopped at $\sigma_h\land\omega$. 

First we prove Side-lemma 5:
$$
\align
&
k
=
\expect\Big(
Y_{\sigma_h\land\omega}
\Big| Y_0=k \Big)
=
\expect\Big(
Y_{\sigma_h}
\Big| [Y_0=k] \,\land\, [\sigma_h<\infty] \Big)
\prob\Big(
\sigma_h<\infty
\Big| Y_0=k \Big)
\\
&
\qquad\qquad\qquad
\le
\left( h + C\sqrt h \right)
\prob\Big(
\sigma_h<\infty
\Big| Y_0=k \Big).
\tag 6.25
\endalign
$$
Where, in the last inequality we applied (6.13). 
Hence (5.23).

\noindent
\qed (Side-lemma 5)
\smallskip

Next we prove Side-lemma 2:
$$
k^2
=
\expect\Big(
Y^2_{\sigma_h\land\omega}
-
2\sum_{s=0}^{\sigma_h\land\omega-1}Y_s
\Big| Y_0=k \Big)
\le
\expect\Big(
Y^2_{\sigma_h\land\omega}
\Big| Y_0=k \Big)
-
2 \expect\Big(
\sigma_h\land\omega
\Big| Y_0=k \Big),
\tag 6.26
$$
where in the last inequality we used the fact that 
$Y_s\ge1$ for $s<\omega$. Hence 
$$
2 \expect\Big(
\sigma_h\land\omega
\Big| Y_0=k \Big)
\le 
\expect\Big(
Y^2_{\sigma_h}
\Big| [Y_0=k] \,\land\, [\sigma_h<\infty] \Big)
\prob\Big(
\sigma_h<\infty
\Big| Y_0=k \Big)
-
k^2
<
C h^2.
\tag 6.27
$$
In the last inequality we used (6.14).
 
\noindent
\qed (Side-lemma 2)
\enddemo

\demo{Proof of Side-lemmas 6 and 7}
\newline
We apply the Optional Stopping Theorem to the martingale 
$Z_t-t$ (for Side-lemma 6),
respectively, to the supermartingale 
$t^2-2tZ_t$ (for Side-lemma 7), 
$t\ge0$, 
both stopped at $\tau_h$.

First Side-lemma 6:
$$
\expect\Big(
\tau_h
\Big| Z_0=k \Big)
=
\expect\Big(
Z_{\tau_h}
\Big| Z_0=k \Big)
-
k
\le
h-k + C\sqrt h.
\tag 6.28
$$
In the last inequality (6.15) was used. 

\noindent
\qed (Side-lemma 6)
\smallskip

Next, Side-lemma 7:
$$
\expect\Big(
\tau_h^2
\Big| Z_0=k \Big)
\le 
2 
\expect\Big(
\tau_hZ_{\tau_h}
\Big| Z_0=k \Big)
\le 
2
\sqrt{
\expect\Big(
\tau_h^2
\Big| Z_0=k \Big)
}
\sqrt{
\expect\Big(
Z_{\tau_h}^2
\Big| Z_0=k \Big)
}
\tag 6.29
$$
Hence, using (6.16) we get
$$
\expect\Big(
\tau_h^2
\Big| Z_0=k \Big)
\le 
4
\expect\Big(
Z_{\tau_h}^2
\Big| Z_0=k \Big)
\le 
C h^2.
\tag 6.30
$$

\noindent
\qed (Side-lemma 7)
\enddemo


\bigskip
\noindent
{\bf Acknowledgement:}
This work was partially supported by the following grants: 
OTKA-T26176 (Hungarian National
Foundation for Scientific Research), 
FKFP-0638/99 (Ministry of Education) and 
TKI-Stochastics\@TUB (Hungarian Academy of Sciences).
\bigskip
\noindent
{\bf References:}
 
\medskip

\item{[1]}
Bass, R.F., Eisenbaum, N., Shi, Z. (2000): 
The most visited sites of symmetric stable process.
{\sl Probability Theory and Related Fields} (to appear)

\item{[2]}
Bass, R.F., Griffin, P.S. (1985):
The most visited site of Brownian motion and random walk,
{\sl Z. Wahrscheinlichkeitstheorie verw. Gebiete} 
{\bf 70}: 417-436

\item{[3]}
Cs\'aki, E., R\'ev\'esz, P., Shi, Z. (2000):
Favourite sites, favourite values and jumping sizes 
for random walk and Brownian motion.
{\sl Bernoulli}, to appear

\item{[4]}
Cs\'aki, E., Shi, Z. (1998):
Large favourite sites of simple random walk and the 
Wiener process.
{\sl  Electronic Journal of Probability} {\bf 3}, 
paper no. 14, pp 1-31.

\item{[5]}
Eisenbaum, N. (1997):
On the most visited sites by a symmetric 
stable process.
{\sl Probability Theory and Related Fieldas} {\bf 107}: 527--535.

\item{[6]}
Erd\H os, P., R\'ev\'esz, P. (1984):
On the favourite points of random walks. 
{\sl Mathematical Structures -- 
Computational Mathematics -- 
Mathematical Modelling.\/} (Sofia)
{\bf 2}: 152-157.

\item{[7]}
Erd\H os, P., R\'ev\'esz, P. (1987):
Problems and results on random walks. 
In: 
{\sl Mathematical Statistics and Probability Theory\/} 
Eds.: P. Bauer, F. Koneczny, W. Wertz.
pp. 59-65. D. Reidel, Dordrecht.

\item{[8]}
Erd\H os, P., R\'ev\'esz, P. (1991): 
Three problems on the random walk in ${\Z}^d$. 
{\sl Studia Sci. Math. Hung.\/} 
{\bf 26}: 309-320.

\item{[9]}
Erd\H os, P. (1994): 
My work with P\'al R\'ev\'esz. 
{\sl Lecture delivered at the
conference dedicated to  the 60th birthday of 
P\'al R\'ev\'esz.}  
Budapest, 8 June 1994.

\item{[10]}
Kesten, H. (1965):
An iterated logarithm law for the local time. 
{\sl Duke Mathematical  Journal} {\bf 32}:  447--456. 

\item{[11]}
Knight, F.B. (1963): 
Random walks and a sojourn density process of 
Brownian motion. 
{\sl Transactions of AMS\/} 
{\bf 109}: 56-86.

\item{[12]}
K\"onig, W. (1996):
A central limit theorem for a one-dimensional polymer measure.
{\sl The Annals of Probability\/} {\bf 24}: 1012-1035.

\item{[13]}
Ray, D. (1963):
Sojourn times of a diffusion process. 
{\sl Illinois Journal of  Mathematics} {\bf 7}: 615--630.

\item{[14]}
R\'ev\'esz, P. (1990): 
{\sl Random Walk in Random and Non-Random Environment.\/} 
World Scientific, Singapore.

\item{[15]}
T\'oth, B (1995):
The true self-avoiding walk with bond repulsion on $\Z$: 
limit thorems.
{\sl The Annals of Probability\/} {\bf 23}: 1523-1556.

\item{[16]}
T\'oth, B (1996):
Generalized Ray-Knight theory and limit theorems for self-interacting random walks on $\Z$.
{\sl The Annals of Probability\/} {\bf 24}: 1324-1367.

\item{[17]}
T\'oth, B., Werner, W. (1997):
Tied favourite edges for simple random walk.
{\sl Combinatorics, Probability and Computing\/}
{\bf 6}: 359-369.

\item{[18]}
Williams, D. (1991):
{\sl Probability with Martingales.}
Cambridge University Press. 
\end